\documentclass{article}
\usepackage[nonatbib]{ml4eng_2020}

\usepackage[utf8]{inputenc} 
\usepackage[T1]{fontenc}    
\usepackage{hyperref}       
\usepackage{url}            
\usepackage{booktabs}       
\usepackage{amsfonts}       
\usepackage{nicefrac}       
\usepackage{microtype}      
\usepackage{amsmath,cite,xspace}
\usepackage{graphicx}
\usepackage{epsfig}
\usepackage{epstopdf}
\usepackage{bm}

\usepackage{color}
\usepackage{amssymb}
\usepackage{hyperref}
\newcommand{\utwi}[1]{\mbox{\boldmath $ #1$}}
\usepackage[table,xcdraw]{xcolor}


\usepackage{mathtools}

\newcommand{\reals}{\mathbb{R}}

\newcommand{\bx}{{\bf x}}

\newcommand{\bv}{{\bf v}}

\newcommand{\bP}{{\bf P}}

\newcommand{\bQ}{{\bf Q}}

\newcommand{\btheta}{{\utwi{\theta}}}

\graphicspath{{./Figures/}} \DeclareGraphicsExtensions{.eps,.pdf,.jpg,.mps,.png}

\title{A Learning-boosted Quasi-Newton Method \\for AC Optimal Power Flow}

 \author{Kyri Baker\\
University of Colorado Boulder\\
Renewable and Sustainable Energy Institute\\
Boulder, CO, USA\\
\texttt{kyri.baker@colorado.edu}
}
\makeatletter
\let\old@ps@headings\ps@headings
\let\old@ps@IEEEtitlepagestyle\ps@IEEEtitlepagestyle
\def\psccfooter#1{%
    \def\ps@headings{%
        \old@ps@headings%
        \def\@oddfoot{\strut\hfill#1\hfill\strut}%
        \def\@evenfoot{\strut\hfill#1\hfill\strut}%
    }%
    \def\ps@IEEEtitlepagestyle{%
        \old@ps@IEEEtitlepagestyle%
        \def\@oddfoot{\strut\hfill#1\hfill\strut}%
        \def\@evenfoot{\strut\hfill#1\hfill\strut}%
    }%
    \ps@headings%
}
\makeatother

\begin{document}
\maketitle

\begin{abstract}
Power grid operators typically solve large-scale, nonconvex optimal power flow (OPF) problems throughout the day to determine optimal setpoints for generators while adhering to physical constraints. Despite being at the heart of many OPF solvers, Newton-Raphson can be slow and numerically unstable. To reduce the computational burden associated with calculating the full Jacobian and its inverse, many Quasi-Newton methods attempt to find a solution to the optimality conditions by leveraging an approximate Jacobian matrix. In this paper, a Quasi-Newton method based on machine learning is presented which performs iterative updates for candidate optimal solutions without having to calculate a Jacobian or approximate Jacobian matrix. The proposed learning-based algorithm utilizes a deep neural network with feedback. With proper choice of weights and activation functions, the model becomes a contraction mapping and convergence can be guaranteed. Results shown for networks up to 1,354 buses indicate the proposed method is capable of finding approximate solutions to AC OPF very quickly.
\end{abstract}

\section{Introduction}

AC optimal power flow (OPF) is a canonical power systems operation problem that is at the heart of optimizing large-scale power networks. Newton-based methods are at the core of many AC OPF solvers used in grids today. However, traditional Newton-Raphson can suffer from slow and numerically unstable Jacobian matrix inversions. To reduce the computational burden associated with calculating the full Jacobian and its inverse, many Quasi-Newton methods attempt to find a solution to the optimality conditions by leveraging an approximate Jacobian matrix  \cite{QN1,QN2}. For example, the full Jacobian can be replaced with a function of the gradient of the objective, or the chord method can be used, which fixes the Jacobian to a constant value from the first iteration. Various ways to construct the Jacobian using data-driven techniques have been explored for power flow in \cite{JacobianData16,JacobianNN} with promising results. Here, we avoid calculating a Jacobian by replacing the Newton-Raphson step with a purely data-driven machine learning (ML) model that learns subsequent iterations. The ML model is trained on Newton-Raphson iterations and learns how to imitate the Newton-Raphson algorithm without having to construct Jacobian matrices or calculate matrix inverses.

Learning for OPF is a rapidly growing area due to the immense power that ML and deep learning in particular can provide for representing extremely complex variable relationships and performing inference (making predictions) extremely quickly. For example, \cite{chatzos2020highfidelity, OPFGNN} attempt to learn AC OPF solutions; \cite{Zamzam_learn_19} develops a method to learn \emph{feasible} AC OPF solutions; ML was used to reduce AC OPF complexity in \cite{hasan2020hybrid,Baker_jointCC}, and warm-start points for AC OPF were obtained using ML in \cite{Baker_learning,Diehl2019WarmStartingAO}, to name a few applications (a more in-depth overview can be found in \cite{Amin20}). Instead of learning the AC OPF solution directly, this paper aims to develop an iterative model that learns Newton-Raphson descent directions depending on the current candidate solution. This technique aims to learn a representation such as ``Given my current step, what direction should I move to decrease the objective function?" rather than ``Given an initial guess, what is the exact optimal solution of the nonconvex AC OPF?" A similar technique was explored in \cite{li2016learning} where a reinforcement learning agent was trained to solve general unconstrained optimization problems with great success. 

The model is comprised of a fully connected three-layer neural network $F_R$ with feedback, where input $\bx^{k}$ is the candidate optimal solution vector at iteration $k$. Reminiscent of a simple recurrent neural network, the model iteratively uses feedback from the output layer as inputs until convergence ($||\bx^{k+1}-\bx^{k}|| \leq \epsilon$). The model thus bypasses any construction of a Jacobian matrix or associated inverse. Simulation results demonstrate that although the solution is approximate, a solution can be obtained orders of magnitudes faster than traditional Newton-Raphson, and numerical issues due to nearly-singular Jacobians are avoided.



\section{Learning-boosted Quasi-Newton} \label{sec:LBQN}
A general nonconvex optimization problem with $n$-dimensional optimization variable vector $\bx$, cost function $f(\cdot): \reals^{n} \rightarrow \reals$, $M$ equality constraints $g_i(x) = 0$, $g_i(\cdot): \reals^{n} \rightarrow \reals$, and $P$ inequality constraints $h_j(x) \leq 0$, $h_j(\cdot): \reals^{n} \rightarrow \reals$ can be written as

\begin{subequations} \label{eqn:gen_opt}
\begin{align} 
\min_{\substack{\bx}} ~~ f(\bx)& \\
\mathrm{s.t:~} g_i(\bx) &= 0, ~~i=1, ... M\\
h_j(\bx) &\leq 0, ~~j=1, ..., P
 \end{align}
\end{subequations}

Newton-Raphson (sometimes just called Newton's Method) utilizes first and second-order derivatives of the Karush Kuhn Tucker (KKT) optimality conditions of \eqref{eqn:gen_opt} to find a stationary point of the KKT conditions. For AC OPF, this results in finding a local (or maybe global) minimum of the nonconvex problem \eqref{eqn:gen_opt}. The Newton-step, that is, the iterative update of candidate optimal solution $\bx^{k+1}$ at iteration $k$ can be written as

\begin{equation}\label{eqn:OGNR}
    \bx^{k+1} = \bx^{k} - \alpha J^{-1}(\bx^{k})d(\bx^{k}),
\end{equation}

\noindent where $d(\bx^{k})$ is a vector of KKT conditions evaluated at the current candidate solution $\bx^{k}$, $J^{-1}(\bx^{k})$ is the Jacobian matrix of the KKT conditions evaluated at $\bx^{k}$, and $\alpha$ is an optional step-size parameter where $0 < \alpha \leq 1$. Equation \eqref{eqn:OGNR} is repeated until convergence, which is typically when iterations cease to change significantly (e.g., $||\bx^{k+1}-\bx^{k}|| \leq \epsilon$) or when the KKT conditions are nearly satisfied (e.g. $d(\bx^{k}) \leq \epsilon$) for some small $\epsilon$, for example.

\subsection{Quasi-Newton Methods}

The term ``Quasi-Newton" simply refers to methods where the Jacobian matrix in \eqref{eqn:OGNR} or its inverse is approximated. These methods are typically used for two reasons: First, calculating the Jacobian and its inverse is expensive and time-consuming, which may not be appropriate for large problems or fast-timescale optimization. Second, the Jacobian can be singular or close to singular at the optimal solution \cite{singularJ05, BakerJacobian13}, making the inverse challenging. Thus, Quasi-Newton methods replace $J^{-1}$ with an approximate Jacobian $\tilde{J}^{-1}$ or sometimes with an approximate Jacobian inverse.

The learning-boosted method replaces  $\alpha^{k}J^{-1}(\bx^{k})d(\bx^{k})$ in \eqref{eqn:OGNR} with a fully connected neural network (NN) model $F_R(\cdot): \reals^n \rightarrow \reals^n$ that takes in $\bx^{k}$ as an input and provides $\bx^{k+1}$ as an output; e.g.

\begin{equation}\label{eqn:LBNR}
    \bx^{k+1} = F_R(\bx^{k}).
\end{equation}

A fully connected three-layer NN is used here. The variable vector $\bx^{k} = [\bv^{k},~ \btheta^k,~ \bP^k_g,~ \bQ^k_g]^T$, where $\bv^{k}$ contains the complex voltage magnitudes at each bus, $\btheta^k$ contains the complex voltage angles, and $\bP_g$ and $\bQ_g$ are the real and reactive power outputs at each generator, respectively. 

\subsection{Network architecture}
Choosing the number of nodes in the hidden layer can be performed by using a popular heuristic \cite{masters93}: $N_h \approx \sqrt{N_i \cdot N_o} = \sqrt{(2N_L+n) \cdot n}$, where $N_i$ is the number of nodes in the input layer and $N_o$ is the number of nodes in the output layer. For the considered problem, the number of inputs is the number of optimization variables $n$ and the number of loads $2N_L$ ($N_L$ real demands and $N_L$ reactive demands). Note that even though the loads are inputs to the NN, they are fixed and do not change as $k$ changes. These numbers were later refined through trial-and-error. A rectified linear unit (ReLU) was used as the activation function on the input layer; a hyperbolic tangent (tanh) activation function was used in the hidden layer, and a thresholded linear function was used on the output.



\section{Convergence Analysis} \label{sec:convergence}

Quasi-Newton iterations are fixed-point iterations. That is, they can be written in the form $\bx^{k+1} = F(\bx^{k})$, where $F(\bx^{k})$ is the right-hand side of \eqref{eqn:OGNR}. Unfortunately, when solving AC OPF it is generally difficult to show $\bx^{k+1} = \bx^{k}$ for any $k$. However, we can design our NN to guarantee convergence.




In \cite{RNNcont}, conditions on the activation functions and bounds on the weights are given to ensure convergence of a recurrent neural network (RNN). The proposed NN architecture is a basic version of an RNN, as the output of the last layer becomes an input of the input layer. Consider the mapping $F_R$ in \eqref{eqn:LBNR}. In order to show that the proposed NN converges to a unique fixed point for any initial $x^0$ and that $F_R$ is a contraction, it must first be true that each activation function $f_i$ must be bounded, continuous, differentiable, real-valued, and have bounded derivatives. In the proposed NN, thresholds on both the input (ReLU) and output (Linear) layers were placed to bound their output. The ReLU function on input $x$, $\max\{0,x\}$, is techncially not differentiable at $x=0$. Thus, a more formal analysis would have to be performed to determine if exactly the same convergence results apply in every situation with a ReLU, although they were achieved in simulation. Similar but slightly worse results were achieved with a tanh function on the input layer, which does satisfy these requirements.

From Theorem 1 and 2 in \cite{RNNcont} we have a bound on each weight $w_{ij}$ connecting nodes $i$ and $j$ that must hold for the network to converge to a unique fixed point for any given $x^0$:

\begin{equation}\label{eqn:weightbound}
    |w_{ij}| < c^* < \frac{1}{N_n \cdot f'_{max}},
\end{equation}

\noindent where $'$ denotes a first-order derivative, $N_n$ is the total number of network nodes, $f_l$ is activation function $l = 1,...,N_n$, and

\begin{equation}
    f_{max} = \max_{\substack{y \in \reals}, l=1,...,N_n}(|f'_l(y)|).
\end{equation}
  
\noindent Then $F_R$ is a contraction with $0 < c < 1$, and $c = N_n c^* |f'_{max}|$, and $w_{ij}$ can be designed by \eqref{eqn:weightbound}. 



\section{Network Architecture} \label{sec:net_params}

Four networks were considered: The IEEE 30 bus, IEEE 300 bus, PG-lib 500-bus \cite{PGlib-OPF}, and 1,354-bus PEGASE networks \cite{pegase}. The number of loads, lines, generators, and total (real) generation capacity for each network is shown in Table \ref{tab:networks}. Line flow constraints were neglected (although some networks did not have them to begin with). The data generation, training and testing of the network, and simulations were performed on a 2017 MacBook Pro laptop with 16 GB of memory. Keras with the Tensorflow backend was used to train the neural network using the Adam optimizer. The chosen number of nodes in the hidden layer and training dataset sizes are shown in Table \ref{tab:hyper}. 

\begin{table}[]\centering
\small \caption{Considered network parameters}
\begin{tabular}{|l|l|l|l|l|}
\hline
\textbf{Case}      & \textbf{\begin{tabular}[c]{@{}l@{}}\# of \\ Loads\end{tabular}} & \textbf{\begin{tabular}[c]{@{}l@{}}\# of \\ Lines\end{tabular}} & \textbf{\begin{tabular}[c]{@{}l@{}}\# of \\ Gens\end{tabular}} & \textbf{\begin{tabular}[c]{@{}l@{}}Total Real \\ Gen. Capacity\end{tabular}} \\ \hline
\textbf{30-bus}    & 20  & 41 & 6 & 335 MW  \\ \hline
\textbf{300-bus}   & 191 & 411 & 69  & 32.68 GW   \\ \hline
\textbf{500-bus}   & 200 & 597  & 56  & 12.19 GW  \\ \hline
\textbf{1,354-bus} & 621 & 1991  & 260  & 128.74 GW  \\ \hline
\end{tabular}\label{tab:networks}
\end{table}

\begin{table}[]\centering\small \caption{Number of nodes and training samples}
\begin{tabular}{|l|l|l|l|l|}
\hline
\textbf{Case}      & \textbf{\begin{tabular}[c]{@{}l@{}}Input\\ Nodes\end{tabular}} & \textbf{\begin{tabular}[c]{@{}l@{}}Output\\ Nodes\end{tabular}} & \textbf{\begin{tabular}[c]{@{}l@{}}Hidden\\ Nodes\end{tabular}} & \textbf{\begin{tabular}[c]{@{}l@{}}Training\\ Samples\end{tabular}} \\ \hline
\textbf{30-bus}    & 112            & 72             & 100   & 72,111    \\ \hline
\textbf{300-bus}   & 1,120          & 738           & 800 & 91,432  \\ \hline
\textbf{500-bus} & 1,512          & 1,112          & 2,300 & 111,674 \\ \hline
\textbf{1,354-bus} & 4,470          & 3,228         & 6,000 & 126,724 \\ \hline
\end{tabular}\label{tab:hyper}
\end{table}

\subsection{Dataset generation}\label{sec:datagen}

The MATPOWER Interior Point Solver (MIPS) \cite{MATPOWER} was used to generate the data and was used as the baseline for comparison with the NN model. The termination tolerance of the MIPS solver was set to $10^{-9}$ for data generation and $10^{-4}$ for testing. The tolerance of the learning-boosted solver was set to $10^{-4}$, where convergence is reached when $||\bx^{k+1}-\bx^k|| \leq \epsilon$. A smaller tolerance was used for data generation to promote smoother convergence and ``basins of attraction" within the ML model. For a fair comparison, the same convergence criteria was used for the NN model and MIPS during testing. 1,000 different loading scenarios were randomly generated at each load bus from a uniform distribution of $\pm40\%$ around the given base loading scenario in MATPOWER. In some cases, the generated load profile resulted in an infeasible solution; these samples were excluded. 

Table \ref{tab:hyper} shows the number of training samples $(\bx^{k}, \bx^{k+1})$ generated for each scenario. Note that a different number of iterations to convergence is encountered for each scenario and network. It is recognized that generating a diverse and representative dataset is an important and essential thrust of research within learning-based OPF methods. This is an important direction of future work.

\section{Simulation Results} \label{sec:sim_result}

\subsection{Prediction results}
1,000 testing scenarios were generated for each network. Each simulation was initialized with a flat start ($\bv^0 = \bP^0_g = \bQ^0_g = 1$ and $\btheta^0 = 0$). 1,000 test scenarios were generated using the same methodology in \ref{sec:datagen}. The mean absolute error (MAE) for voltage magnitudes (in pu) and real power generation (in MW) is shown in Table \ref{tab:error}. The fourth column in the table shows the mean absolute percentage error (MAPE) for the OPF objective function across all scenarios. Interestingly, the 300-bus case had a much worse MAE than other cases. This could be because many solutions in the training dataset had nearly-singular Jacobians near the optimal solution, and thus iterations in some areas were less well-defined and harder for the NN to learn.

\begin{table}[]\centering \small \caption{Mean absolute error across testing dataset}
\begin{tabular}{|l|l|l|l|}
\hline
\textbf{Case}         & \textbf{\begin{tabular}[c]{@{}l@{}}MAE: Voltage\\ Magnitude (pu)\end{tabular}} & \textbf{\begin{tabular}[c]{@{}l@{}}MAE: Active\\ Power (MW)\end{tabular}} & \textbf{\begin{tabular}[c]{@{}l@{}}MAPE: \\Cost (\%)\end{tabular}} \\ \hline
\textbf{30-bus}    & 0.004 pu                                                                                            & 0.64 MW                                                                                       & 0.29\%                                                                                                      \\ \hline
\textbf{300-bus}   & 0.009 pu                                                                                            & 10.47 MW                                                                                      & 0.65\%                                                                                                      \\ \hline
\textbf{500-bus}   & 0.099 pu                                                                                            & 0.62 MW                                                                                       & 0.66\%                                                                                                      \\ \hline
\textbf{1,354-bus} & 0.019 pu                                                                                            & 7.55 MW                                                                                       & 1.16\%                                                                                                      \\ \hline
\end{tabular}\label{tab:error}
\end{table}

\begin{figure}[h!]
    \centering
    \includegraphics[width=0.7\textwidth]{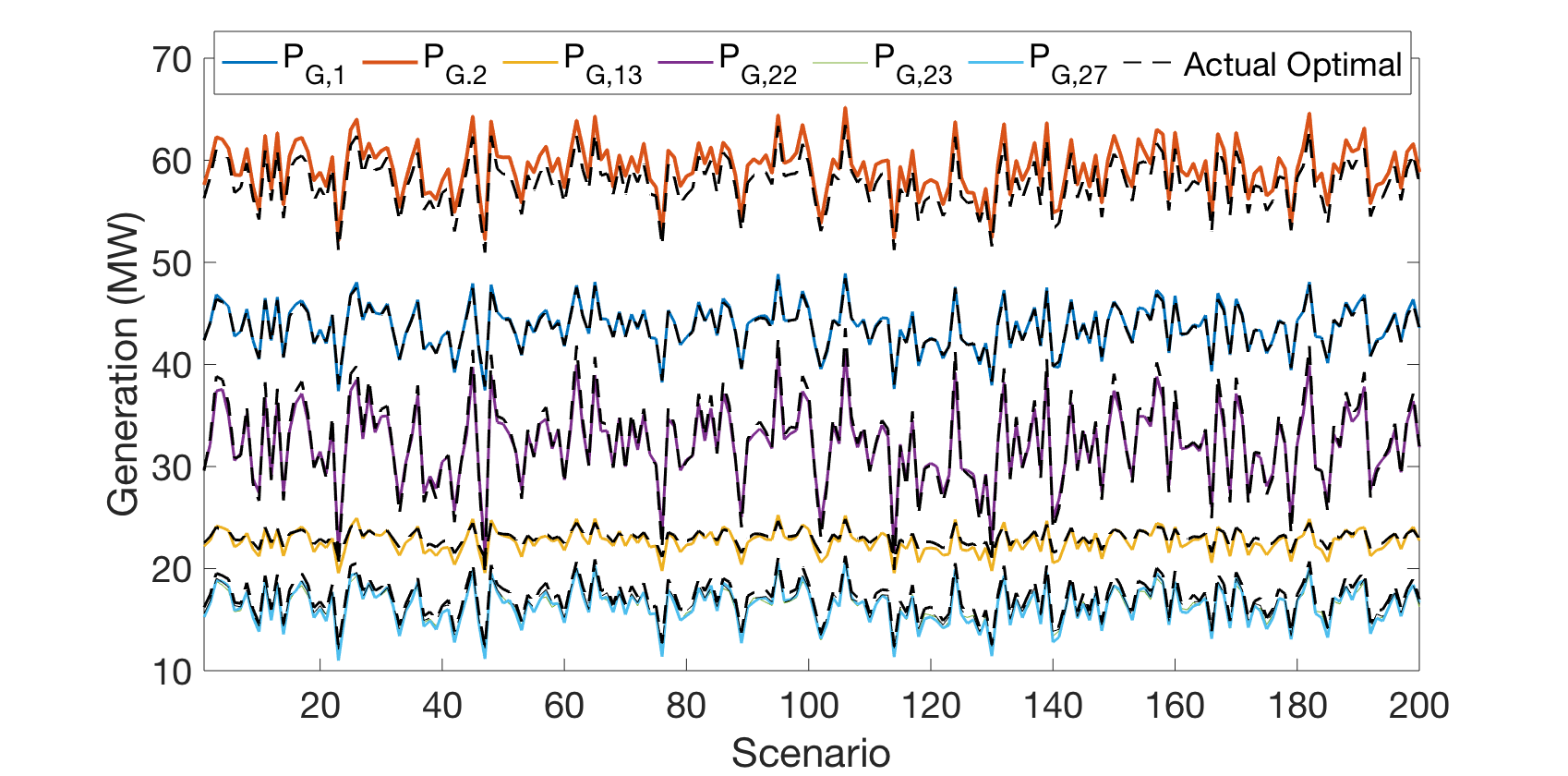}
    \caption{Predicted generation values for 200 test scenarios (colors) and actual optimal values (black dashed line) for the IEEE 30-bus system.}
    \label{fig:gen30}
\end{figure}

Fig. \ref{fig:gen30} shows 200 test scenarios for the 30-bus system. The NN model does a good job of approximating the actual optimal generation values, given as black dashed lines, for each generator. 

\subsection{Computation time}

The mean, worst-case (max), and variance of time to convergence for the MIPS solver and the proposed NN is shown for all networks across all test scenarios in Table \ref{tab:time}.

\begin{table}[h!]\centering \small \caption{Time to convergence for each network.}
\begin{tabular}{|l|l|l|l|l|}
\hline
\textbf{Case}         & \textbf{\begin{tabular}[c]{@{}l@{}}Mean\\ Time (s)\end{tabular}} & \textbf{\begin{tabular}[c]{@{}l@{}}Max\\ Time (s)\end{tabular}} & \textbf{\begin{tabular}[c]{@{}l@{}}Variance \\ in Time (s)\end{tabular}} & \textbf{\begin{tabular}[c]{@{}l@{}}Mean \\ Speedup\end{tabular}} \\ \hline
\textbf{30-bus MIPS}  & 0.04  & 0.41   & 4.52e-04  & \cellcolor[HTML]{656565} \\ \hline
\textbf{30-bus NN}    & 0.06           & 0.42   & 3.53e-04   & \textbf{-0.66x}  \\ \hline
\textbf{300-bus MIPS} & 1.09 & 10.87 & 2.09  & \cellcolor[HTML]{656565}{\color[HTML]{656565} }                  \\ \hline
\textbf{300-bus NN}   & 0.03       & 0.42    & 0.001     & \textbf{36.3x}  \\ \hline
\textbf{500-bus MIPS} & 1.46   & 2.96   & 0.49  & \cellcolor[HTML]{656565}  \\ \hline
\textbf{500-bus NN}   & 0.08       & 0.45      & 3.25e-4  & \textbf{18.3x}   \\ \hline
\textbf{1,354-bus MIPS}   &   7.64     &  19.89     & 15.55 & \cellcolor[HTML]{656565}{\color[HTML]{656565} }     \\ \hline
\textbf{1,354-bus NN}   &   0.34     &  0.69     & 0.0016 & \textbf{22.5x}   \\ \hline
\end{tabular}\label{tab:time}
\end{table}

There is almost no benefit for using a data-driven approach for smaller networks, as these can already be solved in real-time using a solver. However, for larger networks, the learning-boosted Quasi-Newton method obtains solutions extremely quickly compared to the MIPS solver. In addition, the learning-boosted approach has lower variance in computational time, making it more reliable for providing solutions on regular intervals. Interestingly, the 300-bus system takes MIPS longer to solve for some OPF scenarios than the 500-bus. In these cases, the Jacobian was close to singular.

\begin{figure}[h!]
    \centering
    \hspace*{-3mm}\includegraphics[width=0.7\textwidth]{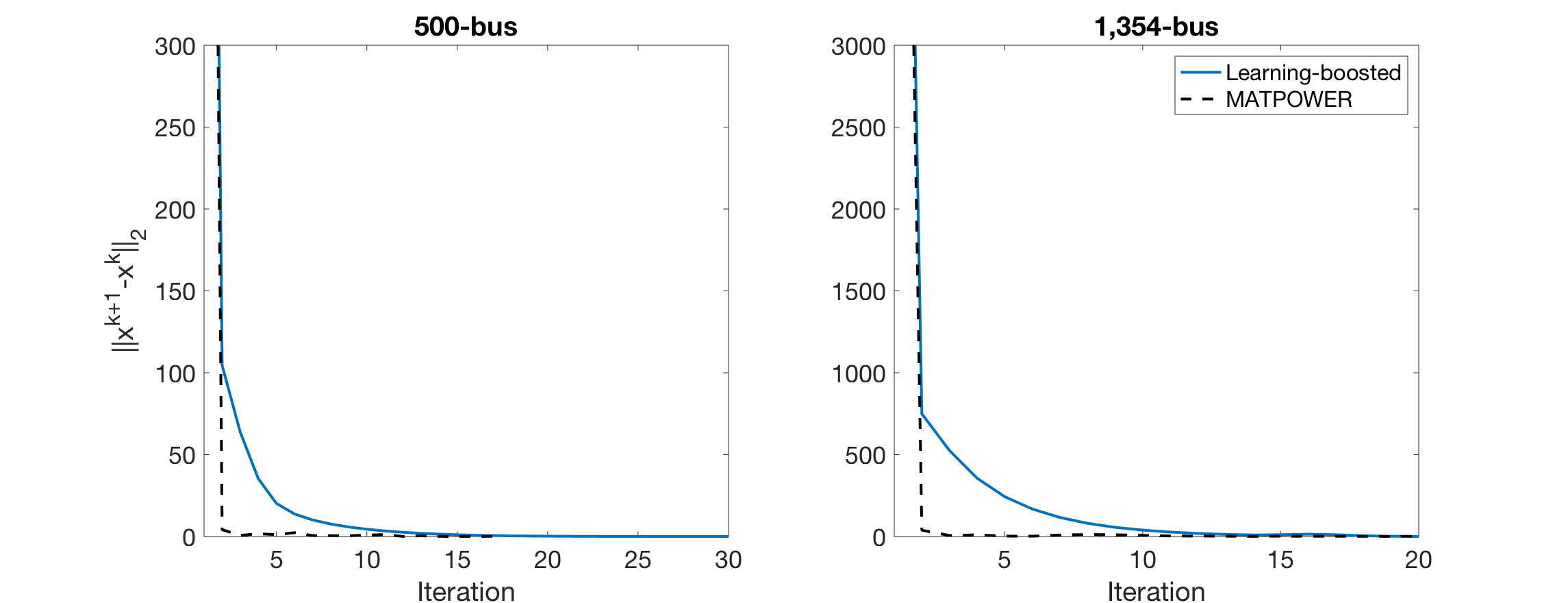}
    \caption{Norm of $\bx^{k+1}-\bx^{k}$ for two scenarios in the 500-bus (left) and 1,354-bus (right) systems.}
    \label{fig:norm_iter}
\end{figure}

Fig. \ref{fig:norm_iter} shows the norm between two iterations, $||\bx^{k+1}-\bx^{k}||_2$, for two typical test scenarios from the 500 and 1,354-bus systems. While MIPS almost always results in fewer iterations to convergence (because it is using the exact Jacobian), each iteration takes longer to perform. Thus, although each iteration in the learning-boosted method is inexact, each iteration is very fast. In addition, convergence of the NN is guaranteed if the conditions are satisfied in Section \ref{sec:convergence}.

\subsection{Assessing feasibility}
Variable bounds are ensured by thresholding the output of the output layer on the NN. However, while convergence of the NN is guaranteed, there is no guarantee of convergence to an AC-feasible solution (just as with DC OPF). In some cases, like with the 30-bus system, the NN actually outputs solutions with smaller mean constraint violations than MIPS. This is likely due to the fact that the training data was generated for a convergence tolerance of $\epsilon = 10^{-9}$ but during testing, the algorithm is terminated at $\epsilon = 10^{-4}$. Despite the prediction error for the 1,354-bus being relatively low (with a mean optimality gap of $1.16\%$ across the testing set), high errors in the satisfaction of the power flow equations indicate a possible need for an increased training dataset size. Table \ref{tab:feasibility} shows the mean constraint violation for the AC power flow constraints across all networks and test scenarios.


\begin{table}[h!]\centering \small \caption{Mean constraint violation across all networks}
\begin{tabular}{|l|l|l|l|l|}
\hline
\textbf{Case}   & \textbf{30-bus} & \textbf{300-bus} & \textbf{500-bus} & \textbf{1,354-bus} \\ \hline
\textbf{\begin{tabular}[c]{@{}l@{}}Mean Constraint \\ Violation (pu)\end{tabular}} & 0.05      & 0.43   &     0.32  & 9.95               \\ \hline
\end{tabular}\label{tab:feasibility}
\end{table}

\subsection{Tracking optimal solutions}
\begin{figure}[t!]
    \centering
    \includegraphics[width=0.65\textwidth]{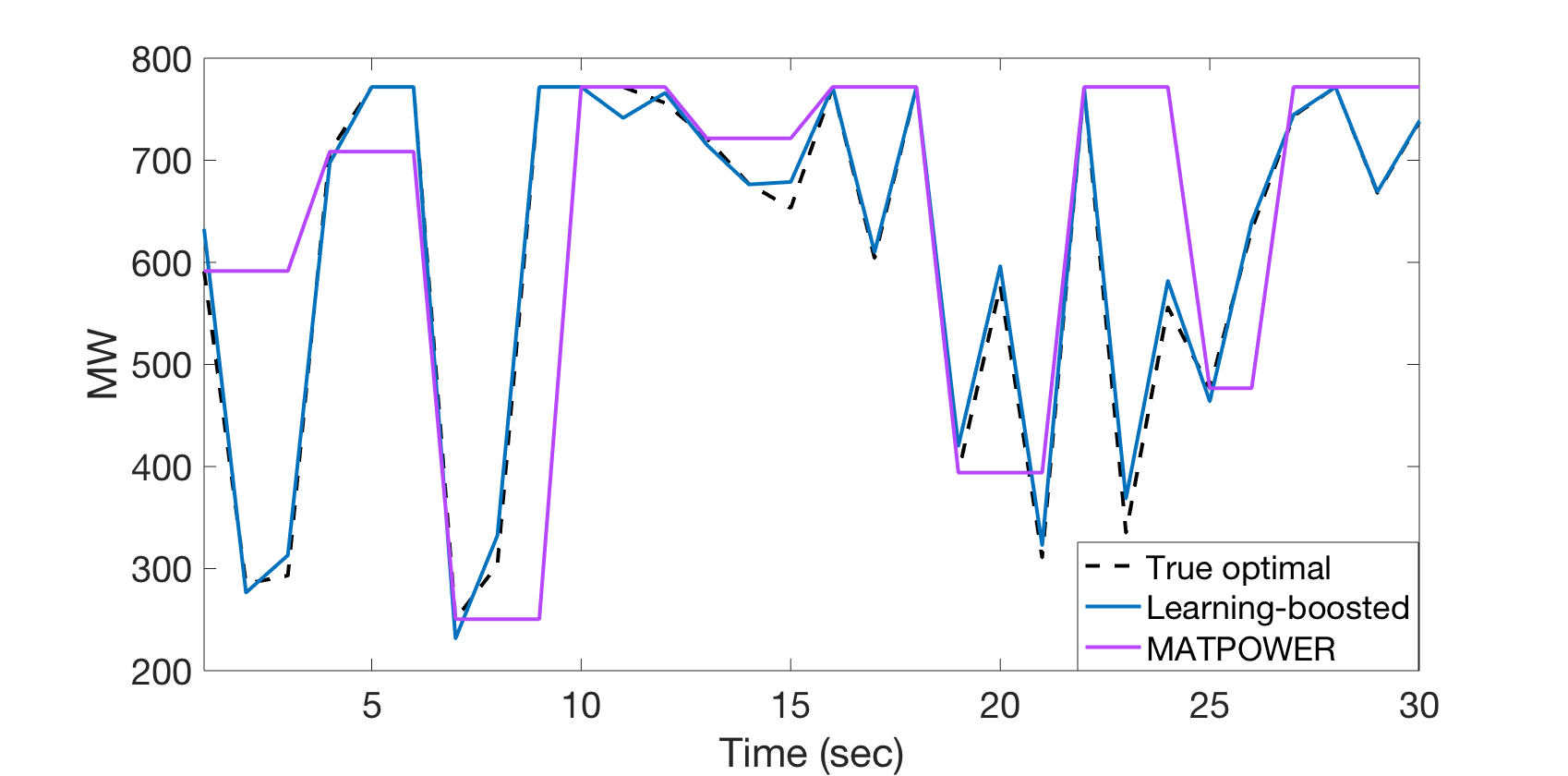}
    \caption{Slack bus generation in the 500-bus system. The learning-boosted solution, although an approximation, is able to track faster changes in demand more accurately than MIPS.}
    \label{fig:500track}
\end{figure}

Assuming measurements of load at each bus are made available on a one-second basis (which is a reasonable assumption; in fact, SCADA systems should provide measurements with latency of less than one second \cite{SCADAdelay}), the benefit of a real-time optimization approach can be further assessed. Typically, real-time adjustments to generators are done via automatic generation control (AGC); a rule-based, suboptimal affine control policy. Here, a simple experiment was performed where both MATPOWER and the learning-boosted algorithm received load updates every one second. Figure \ref{fig:500track} shows the ability of the learning-boosted approach to track optimal generation setpoints with higher accuracy than MATPOWER, despite being an approximation, due to the extremely fast time to obtain new setpoints. This figure illustrates thirty-seconds of the slack bus output for the 500-bus system, which is the generator whose output varies the most in this test case.

As MATPOWER took nearly three seconds to obtain each solution, setpoints were not updated on a one-second scale. Alternatively, the learning-boosted approach provides reasonably accurate solutions in less than a second and keeps up with one-second fluctuations in demand. In reality, generators would adjust their real-time outputs in accordance with the aforementioned AGC heuristics; for illustrative purposes, we show only OPF solutions here. With computationally light methodologies for calculating approximate OPF solutions, optimal solutions could potentially be provided in real-time.

\section{Conclusion} \label{sec:conclusion}

A data-driven model for approximating Newton-Raphson for AC OPF was implemented here. Results show that for smaller networks, not much benefit, other than avoiding singular Jacobian matrices and taking inverses of ill-conditioned matrices, can be found. However, for larger networks, the NN model can provide approximate solutions extremely quickly and at very regular intervals compared to a state-of-the-art OPF solver. The optimality gap was small for all considered cases; however, the feasibility gap for the 1,354-bus system was large. This suggests a need for more training data or an increase in model complexity, as the proposed model only contained one hidden layer. 




\bibliographystyle{IEEEtran}
\bibliography{references.bib}

\end{document}